\documentclass[12pt]{article}
\usepackage[cp1251]{inputenc}
\usepackage[english]{babel}
\usepackage{amsfonts}
\usepackage{amssymb,amsmath,amsthm,eufrak,euscript}
\usepackage{graphicx,fancybox}

\newtheorem{lemma}{Lemma}
\newtheorem{theorem}{Theorem}

\newtheorem{proposition}{Proposition}
\newtheorem{example}{Example}

\begin{document}

\begin{center}
{\bf \Large Polynomial integrals of magnetic geodesic flows on the 2-torus on several energy levels}
\end{center}

\begin{center}
{\bf S.V. Agapov}
\footnote{The author is supported by the grant of the Government of Russia (contract No.~14.Y26.31.0025).}
{\bf, A.A. Valyuzhenich}
\end{center}

\begin{abstract}
In this paper the geodesic flow on a 2-torus in a non-zero magnetic field is considered. Suppose that this flow admits an additional first integral $F$ on $N+2$ different energy levels which is polynomial in momenta of arbitrary degree $N$ with analytic periodic coefficients. It is proved that in this case the magnetic field and metrics are functions of one variable and there exists a linear in momenta first integral on all energy levels.
\end{abstract}

\section{Introduction}

Searching for Riemannian metrics on 2-surfaces with integrable geodesic flows is a classical problem which has been studied extensively for a long time. There are certain topological obstacles to complete integrability, namely, it was proved in~\cite{17} that on any surface of genus $g>1$ there is no analytic Riemannian metrics with integrable geodesic flows. On the other hand, on a 2-sphere and a 2-torus there exist metrics admitting additional polynomial in momenta first integrals.

On the 2-torus which will be considered in this paper there exist at least 2 types of Riemannian metrics with integrable geodesic flow. If metrics has the form $ds^2 = \Lambda (\alpha x + \beta y) (dx^2 + dy^2)$ or $ds^2 = (\Lambda_1 (\alpha_1 x + \beta_1 y) + \Lambda_2 (\alpha_2 x + \beta_2 y)) (dx^2 + dy^2),$ then there exists an additional polynomial in momenta first integral of degree 1 or 2. It is not known whether there exist metrics with irreducible polynomial integrals of higher degree. V.V. Kozlov proposed that on the 2-torus there is no such metrics (see~\cite{11}). In~\cite{DK},~\cite{11} this conjecture was proved under certain assumptions either on the first integral or on the metrics. In general case, however, this conjecture is not proved so far. As it was shown in~\cite{4}, this problem can be reduced to a remarkable quasi-linear system of partial differential equations and the question of existence of smooth global periodic solutions to this system. This system was proved to be semi-Hamiltonian. Such systems have many beautiful properties, they were introduced in~\cite{24} and are of a big interest. The different questions related to such systems and integrable geodesic flows on the 2-torus have been studied in many papers (e.g., see~\cite{4}-~\cite{6},~\cite{16},~\cite{19}).

Geodesic flow in magnetic field on a 2-surface is given by the Hamiltonian system
$$
\dot{x}^j = \{x^j,H\}_{mg}, \qquad \dot{p}_j = \{p_j,H\}_{mg}, \qquad j=1,2 \eqno(1.1)
$$
with Hamiltonian
$H = \frac{1}{2} g^{ij} p_ip_j$
and the Poisson bracket of the following form:
$$
\{F,H\}_{mg} = \sum_{i=1}^2 \left ( \frac{\partial F}{\partial x^i} \frac{\partial H}{\partial p_i} - \frac{\partial F}{\partial p_i} \frac{\partial H}{\partial x^i} \right ) + \Omega (x^1,x^2) \left ( \frac{\partial F}{\partial p_1} \frac{\partial H}{\partial p_2} - \frac{\partial F}{\partial p_2} \frac{\partial H}{\partial p_1} \right ). \eqno(1.2)
$$
If $\{F,H\}_{mg} =0$, then the function $F(x^1,x^2,p_1,p_2)$ is called the first integral of the magnetic geodesic flow (1.1).

The questions related to integrability of geodesic flows in a magnetic field have been studied in many papers (e.g., see~\cite{1},~\cite{3},~\cite{7},~\cite{9},~\cite{10},~\cite{Ef},~\cite{15},~\cite{18}, ~\cite{21},~\cite{23}). The only one known example of Riemannian metric on the 2-torus with integrable magnetic geodesic flow on all energy levels is as follows.

\begin{example}
Suppose that metric has the form
$$ds^2 = \Lambda (y) (dx^2 + dy^2),$$
and magnetic field is $$\omega=-u'(y)dx\wedge dy .$$ Then the magnetic geodesic flow is integrable on all energy levels and the first integral is linear in momenta:
$$F_1 = p_1 + u (y).$$
\end{example}

It was proved in~\cite{1} that if the magnetic geodesic flow on the 2-torus admits a quadratic in momenta first integral on all energy levels, then there exists a linear first integral as in the Example 1. In~\cite{22} it was proved that on the 2-torus there exist no irreducible quadratic in momenta first integrals of magnetic geodesic flow even on two different energy levels.

The main results of this paper are as follows.

\begin{theorem}
Suppose that the geodesic flow (1.1) on the 2-torus in a non-zero magnetic field admits an additional polynomial in momenta first integral $F$ of degree 3 with analytic periodic coefficients on 2 different energy levels $\{H=E_1\}, \{H=E_{2}\}$.
Then the magnetic field and metrics are functions of one variable, there exists a linear in momenta first integral $F_1$ on all energy levels and $F$ can be expressed in terms of $F_1$ and a Hamiltonian.
\end{theorem}

\begin{theorem}
Suppose that the geodesic flow (1.1) on the 2-torus in a non-zero magnetic field admits an additional polynomial in momenta first integral $F$ of degree 4 with analytic periodic coefficients on 3 different energy levels $\{H=E_1\}, \{H=E_{2}\}, \{H=E_{3}\}$.
Then the magnetic field and metrics are functions of one variable, there exists a linear in momenta first integral $F_1$ on all energy levels and $F$ can be expressed in terms of $F_1$ and a Hamiltonian.
\end{theorem}

In general case the following statement holds true.

\begin{theorem}
Suppose that the geodesic flow (1.1) on the 2-torus in a non-zero magnetic field admits an additional polynomial in momenta first integral $F$ of arbitrary degree $N > 4$ with analytic periodic coefficients on $N+2$ different energy levels $\{H=E_1\}, \ldots, \{H=E_{N+2}\}$. Then the magnetic field and metrics are functions of one variable and there exists a linear in momenta first integral $F_1$ on all energy levels.
\end{theorem}

Notice that in the Theorems 1, 2 in addition to existence of the linear first integral $F_1$ on all energy levels we prove that the first integral of degree 3 or 4 can be expressed in terms of the Hamiltonian and $F_1$ (in contrast with the Theorem 3).

This article is organized in the following way. First of all, in section 2 we will prove that polynomial integrability of magnetic geodesic flow on the 2-torus on a finite number of different energy levels (generally speaking, this number depends on a degree of an additional first integral) is equivalent to the integrability on all energy levels. Next, in sections 3, 4 we will prove the Theorems 1, 2. The proof of the general case (the Theorem 3) is given in sections 5, 6.

Due to these results, it is reasonable to consider the question of integrability only on one energy level $\{H=E\}.$ Amazingly, it turns out that on the 2-torus there exist analytic Riemannian metrics such that its geodesic flow in a non-zero magnetic field has an irreducible quadratic in momenta first integral on a fixed energy level (see~\cite{1},~\cite{14}).

Let us also mention about natural mechanical systems which are, due to Maupertuis principle (e.g., see~\cite{8}), in a tight relationship with geodesic flows.
Consider a Hamiltonian system on the 2-torus with Hamiltonian
$$
H = \frac{p_1^2+p_2^2}{2}+V(x,y),
$$
where $V(x,y)$ is a smooth (or analytic) periodic potential. It is known that if $V(x,y) = V_1(y)$ or $V(x,y) = V_1(x)+V_2(y),$ then there exists an additional polynomial in momenta first integral of degree 1 or 2. The question of integrability of this system seems to be more tractable. We refer the reader to~\cite{AA},~\cite{2},~\cite{12},~\cite{13},~\cite{20}.

The authors thank Andrey Mironov for interest to this work and very helpful discussions.

\section{Integrability on several energy levels}

First of all let us prove that integrability of magnetic geodesic flow on two different energy levels via cubic in momenta first integral implies the integrability of this flow on all energy levels. Then we will show that via similar arguments (with minimal changes) it can be proved that existence of the first integral of an arbitrary degree $N$ on $N+2$ different energy levels also implies the integrability on all energy levels.

Throughout this paper we will be working in the conformal coordinates $(x,y),$ in which
$ds^2 = \Lambda (x,y) (dx^2+dy^2)$, $H = \frac{p_1^2+p_2^2}{2 \Lambda}.$ As known, it is possible to define these coordinates on the 2-torus globally.

The following statement holds true.

\begin{lemma}
Suppose that the geodesic flow (1.1) on the 2-torus in a non-zero magnetic field admits an additional polynomial in momenta first integral $F$ of degree 3 on two different energy levels $\{H=E_1\}$ and $\{H=E_2\}$ where $E_1 \neq E_2$. Then $F$ is the first integral of the same flow on all energy levels.
\end{lemma}

Let us prove this lemma.
On a fixed energy level $\{H=\frac{C}{2}\}$ let us parameterize the momenta in the following way: $p_1 = \sqrt{C \Lambda} \cos \varphi, p_2 = \sqrt{C \Lambda} \sin \varphi.$
The equations (1.1) take the form
$$
\dot{x} = \sqrt{\frac{C}{\Lambda}} \cos \varphi,
\qquad \dot{y} = \sqrt{\frac{C}{\Lambda}} \sin \varphi,
\qquad \dot{\varphi} = \frac{\sqrt{C}}{2 \Lambda \sqrt{\Lambda}} (\Lambda_y \cos \varphi - \Lambda_x \sin \varphi) - \frac{\Omega}{\Lambda}.
$$
We shall search for the first integral $F$ of the third degree in the following form:
$$F(x,y,\varphi) =a_{0}(x,y)p_1^3+a_{1}(x,y)p_1^2p_2+a_2(x,y)p_1p_2^2+a_3(x,y)p_2^3+$$
$$+b_0(x,y)p_1^2+b_1(x,y)p_1p_2+b_2(x,y)p_2^2+c_0(x,y)p_1+c_1(x,y)p_2+d_0(x,y), \eqno(2.1)$$
where $p_1 = \sqrt{C \Lambda} \cos \varphi, p_2 = \sqrt{C \Lambda} \sin \varphi.$
The condition $\dot{F}=0$ is equivalent to the following
equation:
$$
F_x \cos \varphi + F_y \sin \varphi +
F_{\varphi} \left( \frac{\Lambda_y}{2 \Lambda} \cos \varphi -
\frac{\Lambda_x}{2 \Lambda} \sin \varphi - \frac{\Omega}{\sqrt{C \Lambda}} \right)
= 0. \eqno(2.2)
$$
Let us substitute (2.1) into (2.2) and equate the coefficients at
$e^{i k \varphi}$ to zero. For $k=4$ we obtain:
$$
(a_0-a_2)_x - (a_1-a_3)_y - i ((a_0-a_2)_y +(a_1-a_3)_x) = 0.
$$
Since all the coefficients of $F$ are real functions, we obtain
$$
(a_0-a_2)_x = (a_1-a_3)_y, \qquad (a_0-a_2)_y = -(a_1-a_3)_y,
$$
which implies that $\triangle(a_2-a_0) = \triangle(a_3-a_1) = 0,$
where $\triangle = \partial_{xx} + \partial_{yy}$ is the Laplace operator. Due to the maximum principle for harmonic functions, it follows from here that:
$$
a_2-a_0 = A_0, \ a_3-a_1 = A_1,
$$
where $A_0, A_1$ are arbitrary constants. Applying the appropriate
rotation in the plane $(x,y)$ and dividing $F$ by the appropriate
constant we can assume that $A_0 = 1, A_1 = 0.$

Let us introduce new functions: $f(x,y) = b_0 - b_2, \ g(x,y) = b_1.$
For $k=3$ we have:
$$3 i \Omega - i f_y - g_y + f_x - i g_x = 0$$
which implies that
$$\Omega(x,y) = \frac{1}{3}(f_y+g_x), \qquad f_x = g_y. \eqno(2.3)$$

For $k=2$ we have:
$$6 {c_0}_y - 6 i {c_1}_y + 6 C \Lambda \left({a_0}_y - i {a_1}_y + i {a_0}_x + {a_1}_x\right) + 6 \left(C (1+a_0 - i a_1) \Lambda_y + i {c_0}_x + {c_1}_x\right) + $$
$$+4 (f-i g) \left(f_y + g_x\right) + 3 C (-i+2 i a_0 + 2 a_1) \Lambda_x = 0$$
which implies that
$$3 {c_0}_y + 3 {c_1}_x + 2 f \left(f_y + g_x\right) + 3 C \left((1+a_0) \Lambda_y + \Lambda \left({a_0}_y + {a_1}_x\right) + a_1 \Lambda_x\right) = 0, \eqno(2.4)$$
$$-6 {c_1}_y - 6 C a_1 \Lambda_y + 6 C \Lambda \left(-{a_1}_y + {a_0}_x\right) + 6 {c_0}_x - 4 g \left(f_y+g_x\right) + 3 C (-1+2 a_0) \Lambda_x = 0. \eqno(2.5)$$
Proceeding in this way, for $k=1$ we obtain:
$$
4 c_1 \left(f_y+g_x\right)+4 C a_1 \Lambda \left(f_y+g_x\right)-
$$
$$-3 \left(2 C g \Lambda_y+C \Lambda \left(g_y+4 {b_2}_x+3 f_x\right)+4 \left({d_0}_x+C ({b_2}+f) \Lambda_x\right)\right) = 0, \eqno(2.6)
$$
$$
6 {d_0}_y+2 c_0 \left(f_y+g_x\right)+2 C \Lambda \left(3 {b_2}_y+(1+a_0) \left(f_y+g_x\right)\right)+3 C \left(2 b_2 \Lambda_y+g \Lambda_x\right) =0. \eqno(2.7)
$$
For $k=0$ we have:
$$
\Lambda \left(2 \left({c_1}_y+2 C a_1 \Lambda_y+C \Lambda
\left({a_1}_y+{a_0}_x\right)+{c_0}_x\right)+C (1+4 a_0) \Lambda_x\right)+2 \left(c_1 \Lambda_y+c_0 \Lambda_x\right)=0. \eqno(2.8)
$$

The equations (2.3) --- (2.8) give the necessary and sufficient conditions on $F$ to be the first integral of magnetic geodesic flow on a fixed energy level $H = \frac{C}{2}.$ Notice that the equations (2.4) --- (2.8) are linear in $C.$ Since we require $F$ to be the first integral on two different energy levels $\{C = 2 E_1\}, \{C = 2 E_2\},$ the equations (2.4) --- (2.8) imply that
$$
3 {c_0}_y+3 {c_1}_x+2 f \left(f_y+g_x\right)=0,
$$
$$
(1+a_0) \Lambda_y + \Lambda \left({a_0}_y + {a_1}_x\right)+a_1 \Lambda_x=0,
$$
$$
-3 {c_1}_y+3 {c_0}_x-2 g \left(f_y+g_x\right)=0,
$$
$$
-2 a_1 \Lambda_y+2 \Lambda \left(-{a_1}_y+{a_0}_x\right)+(-1+2 a_0) \Lambda_x=0,
$$
$$
-3 {d_0}_x+c_1 \left(f_y+g_x\right)=0,
$$
$$
-6 g \Lambda_y+\Lambda \left(-3 \left(g_y+4 {b_2}_x+3 f_x\right)+4 a_1 \left(f_y+g_x\right)\right)-12 (b_2+f) \Lambda_x=0,
$$
$$
3 {d_0}_y+c_0 \left(f_y+g_x\right) = 0,
$$
$$
6 b_2 \Lambda_y+2 \Lambda \left(3 {b_2}_y+(1+a_0) \left(f_y+g_x\right)\right)+3 g \Lambda_x =0,
$$
$$
c_1 \Lambda_y+\Lambda \left({c_1}_y+{c_0}_x\right)+c_0 \Lambda_x = 0,
$$
$$
4 a_1 \Lambda_y+2 \Lambda \left({a_1}_y+{a_0}_x\right)+(1+4 a_0) \Lambda_x =0.
$$
Notice that if at the same time $F$ is the first integral on any additional energy level $\{H=E_3\},$ where $E_2 \neq E_3 \neq E_1,$ then we obtain the same system of equations. It means that for cubic in momenta first integral $F$ the integrability on two different energy levels is equivalent to the one on all energy levels. In addition, one can check straightforwardly that the last system of equations is equivalent to the equations (3.6) --- (3.8), (3.11) --- (3.19) from the next section which deals with cubic first integral on all energy levels.
Lemma 1 is proved.

In general case the following statement holds true.

\begin{lemma}
Suppose that the geodesic flow (1.1) on the 2-torus in a non-zero magnetic field admits an additional polynomial in momenta first integral $F$ of arbitrary degree $N$ on $N+2$ different energy levels $\{H=E_1\}, \ldots, \{H=E_{N+2}\}$. Then $F$ is the first integral of the same flow on all energy levels.
\end{lemma}

The idea of the proof of lemma 2 is as follows.
Notice that in the case of polynomial first integral of an arbitrary degree $N$ on different energy levels we obtain the system of equations polynomial in $\sqrt{C}$ which is analogous to the system (2.4) --- (2.8). The maximal possible degree (in $\sqrt{C}$) of equations of this system is $N+1.$ Thus, if one requires the integrability on $N+2$ different energy levels, then all the coefficients of each equation must vanish identically. This implies the equivalence of polynomial integrability on $N+2$ and all energy levels.

\

{\bf Remark. }
We would like to notice that the estimate of a number of energy levels in lemma 2 is rough.
For example, in the case of $N=3$ it is sufficient to require the integrability only on 2 different energy levels (lemma 1). In the case of $N=4$ it is easy to show that 3 different energy levels are enough. In general case, for a first integral $F$ of arbitrary degree $N,$ depending on the parity,  $\frac{N+1}{2}$ or $\frac{N+2}{2}$ energy levels seem to be enough for integrability on all energy levels.

\section{Integral of degree 3}

In this section we give the proof of the Theorem 1.
Suppose that there exists an additional cubic in momenta first integral on all energy levels:
$$F_3=a_{0}(x,y)p_1^3+a_{1}(x,y)p_1^2p_2+a_2(x,y)p_1p_2^2+a_3(x,y)p_2^3+$$
$$+b_0(x,y)p_1^2+b_1(x,y)p_1p_2+b_2(x,y)p_2^2+c_0(x,y)p_1+c_1(x,y)p_2+d_0(x,y),$$
all the coefficients are supposed to be analytic periodic functions in both variables $x,y.$ Magnetic Poisson bracket $\{F_3,H\}_{mg}$ of the form (1.2) is a polynomial in momenta of the fourth degree. Equating coefficients of this polynomial at different degrees to zero we obtain the system of PDEs on metrics $\Lambda(x,y)$ and coefficients of the first integral $F_3.$

The first group of the equations (3.1) --- (3.5) is the same as for standard geodesic flow with no magnetic field:
$$
a_1\Lambda_y+3a_0\Lambda_x+2\Lambda {a_0}_x=0, \eqno(3.1)
$$
$$
a_2\Lambda_y+a_1\Lambda_x+\Lambda ({a_0}_y+{a_1}_x)=0, \eqno(3.2)
$$
$$
(a_1+3a_3)\Lambda_y+(3a_0+a_2)\Lambda_x+2\Lambda ({a_1}_y+{a_2}_x)=0, \eqno(3.3)
$$
$$
a_2\Lambda_y+a_1\Lambda_x+\Lambda ({a_2}_y+{a_3}_x)=0, \eqno(3.4)
$$
$$
3a_3\Lambda_y + a_2\Lambda_x + 2\Lambda {a_3}_y = 0. \eqno(3.5)
$$
Acting by the usual way (following~\cite{16}), let us take $(3.1)-(3.3)+(3.5)$ and $(3.2)-(3.4)$ and by the cross differentiation we obtain the following relations:
$$
\triangle(a_2-a_0) = \triangle(a_3-a_1) = 0,
$$
where $\triangle = \partial_{xx} + \partial_{yy}$ is the Laplace operator. The last relation implies
$$
a_2-a_0 = A_0, \ a_3-a_1 = A_1,
$$
where $A_0, A_1$ are arbitrary constants. Without loss of generality we may assume that $A_0 = 1, A_1 = 0.$
Then the equations (3.1) --- (3.5) take the form
$$
a_1\Lambda_y+2\Lambda {a_0}_x+3a_0\Lambda_x=0, \eqno(3.6)
$$
$$
(1+a_0)\Lambda_y+\Lambda ({a_0}_y+{a_1}_x)+a_1\Lambda_x=0, \eqno(3.7)
$$
$$
3a_1\Lambda_y + 2\Lambda {a_1}_y + (1+a_0)\Lambda_x=0. \eqno(3.8)
$$
Let us denote
$$f(x,y) = b_0(x,y) - b_2(x,y), \qquad g(x,y) = b_1(x,y).$$
Then the second group of the equations takes the following form:
$$
g\Lambda_y+2(b_2+f)\Lambda_x+2\Lambda({b_2}_x+f_x-a_1 \Omega)=0, \eqno(3.9)
$$
$$
g\Lambda_x+2b_2\Lambda_y+2\Lambda(g_x+f_y+{b_2}_y+(a_0-2) \Omega)=0, \eqno(3.10)
$$
$$
g\Lambda_y+2(b_2+f)\Lambda_x+2\Lambda({b_2}_x+g_y-a_1 \Omega)=0, \eqno(3.11)
$$
$$
g\Lambda_x+2b_2\Lambda_y+2\Lambda({b_2}_y+(a_0+1) \Omega)=0, \eqno(3.12)
$$
Taking $(3.10)-(3.12)$ and $(3.9)-(3.11)$ we obtain:
$$
\Omega(x,y) = \frac{1}{3}\left(f_y+g_x\right), \eqno(3.13)
$$
$$
f_x = g_y. \eqno(3.14)
$$
Equations of the third group are as follows:
$$
c_0\Lambda_x+2\Lambda{c_0}_x+c_1\Lambda_y-\frac{2}{3}g\Lambda(f_y+g_x)=0,
$$
$$
3{c_0}_y+3{c_1}_x+2f(f_y+g_x)=0,
$$
$$
c_0\Lambda_x+2\Lambda{c_1}_y+c_1\Lambda_y+\frac{2}{3}g\Lambda(f_y+g_x)=0.
$$
They can be written in another form:
$$
3{c_0}_x-3{c_1}_y-2g(f_y+g_x)=0, \eqno(3.15)
$$
$$
3{c_0}_y+3{c_1}_x+2f(f_y+g_x)=0, \eqno(3.16)
$$
$$
(c_0\Lambda)_x+(c_1\Lambda)_y=0, \eqno(3.17)
$$
Finally, the last group of equations is as follows:
$$
3 {d_0}_x - c_1 (f_y+g_x) = 0, \eqno(3.18)
$$
$$
3 {d_0}_y + c_0 (f_y+g_x) = 0. \eqno(3.19)
$$

The equations (3.6)---(3.8), (3.11)---(3.19) give the necessary and sufficient conditions for $F_3$ to be the first integral of magnetic geodesic flow (1.1).

Consider the equations (3.15), (3.16). By cross differentiation we obtain:
$$
3 \triangle c_0 + (2f(f_y+g_x))_y-(2g(f_y+g_x))_x = 0,
$$
$$
3 \triangle c_1 + (2f(f_y+g_x))_x+(2g(f_y+g_x))_y = 0.
$$
Due to (3.14), we have
$$
f(f_y+g_x) = \left(\frac{f^2}{2}\right)_y + (fg)_x - f_xg = \left(\frac{f^2-g^2}{2}\right)_y + (fg)_x,
$$
$$
g(f_y+g_x) = (fg)_y - fg_y + \left(\frac{g^2}{2}\right)_x  = -\left(\frac{f^2-g^2}{2}\right)_x + (fg)_y.
$$
Taking this into account, we may conclude that (3.15), (3.16) are equivalent to:
$$
\triangle (3 c_0 +  f^2 - g^2) = 0, \quad \triangle (3 c_1 + 2 fg) = 0,
$$
and it follows from here that
$$
3 c_0 +  f^2 - g^2 = K_1, \qquad 3 c_1 + 2 fg = K_2, \eqno(3.20)
$$
where $K_1, K_2$ are arbitrary constants.

Consider the equations (3.18), (3.19). By cross differentiation we obtain:
$$
\left(c_0 (f_y+g_x)\right)_x + \left(c_1 (f_y+g_x)\right)_y = 0. \eqno(3.21)
$$
Taking (3.14), (3.20) into account, we may rewrite (3.21) in the following way:
$$
\left((K_1 + g^2 - f^2) (f_y+g_x)\right)_x + \left((K_2 - 2 fg) (f_y+g_x)\right)_y = 0.
$$
$$
\left((K_1 + g^2 - f^2) (f_y+g_x) + 2 fg (g_y - f_x) \right)_x +
$$
$$
+ \left((K_2 - 2 fg) (f_y+g_x) + (f^2 - g^2) (f_x - g_y)\right)_y = 0.
$$
$$
\left(K_1 (f_y+g_x) + (fg^2)_y - (gf^2)_x + \frac{1}{3} (g^3)_x - \frac{1}{3} (f^3)_y \right)_x +
$$
$$
+ \left(K_2 (f_y+g_x) - (gf^2)_y - (fg^2)_x + \frac{1}{3} (g^3)_y + \frac{1}{3} (f^3)_x \right)_y = 0.
$$
The last relation is equivalent to $\triangle \left(K_1 g +K_2 f +\frac{1}{3} g^3 - g f^2 \right) = 0.$
Hence,
$$
K_1 g +K_2 f +\frac{1}{3} g^3 - g f^2 = K_3 \eqno(3.22)
$$
with a constant $K_3.$
Now we may rewrite (3.14) in the following way:
$$
g_y - f'(g) g_x = 0,
$$
here $f'(g)$ can be found explicitly from (3.22), namely, we have:
$$
f(g) = \frac{3 K_2 \pm \sqrt{9 K_2^2 + 12 g \left(-3 K_3 + 3 K_1 g + g^3\right)}}{6 g},
$$
$$
f'(g) = \frac{ 3 K_2 \left(\pm \sqrt{3} K_2 - Q \right) \mp 2 \sqrt{3} g \left(2 g^3+3 K_3\right)}{6 g^2 Q},
$$
where $Q = \sqrt{3 K_2^2+4 g \left(g^3+3 g K_1 - 3 K_3\right)}.$ Denote $$\alpha(g) = 3 K_2 \left(\pm \sqrt{3} K_2 - Q \right) \mp 2 \sqrt{3} g \left(2 g^3+3 K_3\right), \qquad \beta(g) = 6 g^2 Q,$$
so (3.14) can be rewritten in the following way:
$$
\beta(g)g_y - \alpha(g) g_x = 0. \eqno(3.23)
$$
Let us show that similarly to Hopf equation, the equation (3.23) does not admit non-constant global analytic periodic solutions. Characteristics of (3.23) are straight lines on the plane $(x,y)$ having a velocity vector $(-\alpha(g),\beta(g))$ and given by the equation $$\beta(g) x + \alpha(g)y = c,$$ where $c \in \mathbb{R}$ is an arbitrary constant. Take a smooth curve $\gamma(t) = (x(t),y(t))$ which intersects the characteristics transversally. Consider a Cauchy problem for (3.23) with the initial data $$g\mid_{\gamma(t)} = \phi(x,y)$$ with an arbitrary function $\phi(x,y).$ Suppose that $\phi(x,y)$ is not a constant. Let us choose two points $\xi \neq \eta$ on the curve $\gamma$ such that $0 \neq \phi(\xi) \neq \phi(\eta) \neq 0.$ Then characteristics going from these two points intersect and it means that a solution blows up. Thus $g\mid_{\gamma(t)} = \phi(x,y)$ is a constant and it means that $g(x,y) \equiv G$ is a constant function.

Then, due to (3.14) we have $f(x,y) = f_1(y),$ hence the magnetic field $\Omega$ depends only on one variable:
$$\Omega (y) = \frac{1}{3}f_1'(y).$$
It follows from (3.20) that $c_0(x,y), \ c_1(x,y)$ also depend only on $y:$
$$
c_0 (y) = \frac{1}{3} \left( K_1 + G^2 - f_1(y)^2 \right), \quad
c_1 (y) = \frac{1}{3} \left( K_2 - 2 G f_1(y) \right).
$$
Due to this, (3.18) takes the form $3 {d_0}_x - \frac{f_1'(y)}{3}\left( K_2 - 2 G f_1(y) \right) = 0.$
All the unknown functions are periodic, hence $d_0(x,y) = d_0(y)$ and
$$
f_1'(y) \left( K_2 - 2 G f_1(y) \right) = 0.
$$
This implies that
$$
\frac{\partial}{\partial y}\left( K_2 f_1(y) - G f_1^2(y) \right) = 0,
$$
consequently $K_2 f_1(y) - G f_1^2(y) = C_0$ with an arbitrary constant $C_0.$ If $K_2 \neq 0$ or $G \neq 0,$ then the only solutions of this equation are constants. But in this case the magnetic field $\Omega$ vanishes identically and we have a standard geodesic flow.
Thus it is left to consider the case $K_2 = G = 0.$

\

{\bf Remark. }
Another way to obtain the relations $K_2 = G = 0$ is as follows. Suppose that $f_1(y)$ is not a constant function. Then, due to (3.22), we obtain $K_2 = G = K_3 = 0.$

\

It follows from (3.19) that $d_0(y) = \frac{1}{27} \left( f_1(y)^3 - 3 K_1 f_1(y) \right).$

Thus we have $K_2 = K_3 = G =0.$
The equation (3.17) implies that
$$\left(K_1 - f_1(y)^2\right) \Lambda_x = 0. \eqno(3.24)$$
Suppose that there exists a point $q=(x_0,y_0)$ such that $\Lambda_x \neq 0$ at $q.$ It means that there is a neighborhood $U_q$ of $q$ such that $\Lambda_x \neq 0$ in $U_q.$ Then, due to (3.24) $f_1(y) \equiv \sqrt{K_1}$ or $f_1(y) \equiv -\sqrt{K_1}$ in $U_q.$ If an analytic function $f_1(y)$ is a constant on an open subset, then it is identically constant on the whole domain. But in this case the magnetic field $\Omega$ vanishes identically. So we obtain $\Lambda_x \equiv 0$ from what follows that the metrics depends only on one variable:
$$\Lambda(x,y) = \Lambda(y).$$

Let us show that in this case the first integral $F_3$ can be expressed in terms of Hamiltonian $H$ and a linear in momenta first integral on all energy levels as in the Example 1. The first group of the equations (3.6) --- (3.8) takes the following form:
$$
a_1\Lambda_y+2\Lambda {a_0}_x=0, \eqno(3.25)
$$
$$
(1+a_0)\Lambda_y+\Lambda ({a_0}_y+{a_1}_x)=0, \eqno(3.26)
$$
$$
2\Lambda {a_1}_y + 3a_1\Lambda_y=0. \eqno(3.27)
$$
It follows from (3.25), (3.27) that $3 {a_0}_x = {a_1}_y.$
Due to this, (3.25) can be written as follows
$$a_1\Lambda_y+2\Lambda {a_0}_x = (a_1\Lambda)_y - {a_1}_y\Lambda + 2\Lambda {a_0}_x =$$
$$(a_1\Lambda)_y - \Lambda {a_0}_x = (a_1\Lambda)_y - \Lambda (1+a_0)_x = (a_1 \Lambda)_y - ((1+a_0) \Lambda)_x = 0.$$
Equation (3.26) can be written in the following way:
$$(a_1 \Lambda)_x + ((1+a_0) \Lambda)_y = 0.$$
The last two equalities imply that $\triangle (a_1 \Lambda)=0$ and $\triangle ((1+a_0) \Lambda) = 0,$
Hence,
$$a_1 \Lambda=s_1, \quad (1+a_0) \Lambda = s_0,$$
here $s_0, s_1$ are some constants. So we have
$a_1(y) = \frac{s_1}{\Lambda(y)}, \ a_0(y) = \frac{s_0}{\Lambda(y)}-1.$
Due to $3 {a_0}_x = {a_1}_y$ we have ${a_0}_x={a_1}_y=0,$ hence $a_1$ is a constant. Due to (3.25) we obtain $a_1 = 0.$

The equations (3.11), (3.12) imply that
$$\Lambda(y) {b_2}_x=0, \quad \frac{2}{3} s_0 f_1'(y) + 2(b_2 \Lambda)_y = 0.$$
It follows from the first equation that $b_2(x,y) = b_2(y);$ then the second equation can be integrated:
$b_2(y) = \frac{3 s_2 - s_0 f_1(y)}{3 \Lambda_1(y)},$
here $s_2$ is a constant.

Thus all the coefficients of the first integral $F_3$ are found; it commutes with the Hamiltonian and has the following form:
$$F_3 = \frac{(p_1^2+p_2^2)(3 (p_1s_0+s_2)-s_0 f_1(y))}{3 \Lambda(y)} - $$
$$-\frac{1}{27}(3 p_1 - f_1(y)) \left( 9 p_1^2 - 6 p_1 f_1(y) + f_1(y)^2 - 3 K_1 \right).$$
Recall that the magnetic field is $\Omega(y) = \frac{f_1'(y)}{3}.$
As known, if metrics and magnetic field depend only on one variable, then there exists a linear first integral $F_1 = p_1 - \frac{f_1(y)}{3}$ on all energy levels.
By straightforward calculations one may verify that
$$F_3 = -F_1^3 + 2 s_0 H F_1 + 2 s_2 H + \frac{K_1}{3} F_1,$$
and it means that $F_3$ can be expressed in terms of the Hamiltonian and the linear first integral $F_1.$
Combining this result with the statement of lemma 1, we obtain the proof of Theorem 1.

\section{Integral of degree 4}

The case of an integral of fourth degree is equivalent to the one of third degree but technically it is more complicated. We will skip the bulky calculations and give only the main steps of the proof.

Suppose that there exists an additional polynomial in momenta first integral of fourth degree on all energy levels:
$$F_4=a_{0}(x,y)p_1^4+a_{1}(x,y)p_1^3p_2+a_2(x,y)p_1^2p_2^2+a_3(x,y)p_1p_2^3+a_{4}p_2^4+$$
$$+b_{0}(x,y)p_1^3+b_{1}(x,y)p_1^2p_2+b_2(x,y)p_1p_2^2+b_3(x,y)p_2^3+$$
$$+c_0(x,y)p_1^2+c_1(x,y)p_1p_2+c_2(x,y)p_2^2+d_0(x,y)p_1+d_1(x,y)p_2+e_0(x,y),$$
all the coefficients are supposed to be analytic periodic functions in both variables $x,y.$
Let us fix the Kolokol'tsov constants:
$$
a_2-a_0-a_4 = A_0 = 1, \ a_3-a_1 = A_1 = 0.
$$
and make the change of variables:
$$f(x,y) = b_0 - b_2, \ g(x,y) = b_1 - b_3.$$
Magnetic field has the form: $$\Omega (x, y) = \frac{1}{4} (f_y+g_x).$$
The equations of the first group are as follows:
$$
a_1\Lambda_y+2\Lambda {a_0}_x+4a_0\Lambda_x=0, \eqno(4.1)
$$
$$
2(1+a_0+a_4)\Lambda_y+2 \Lambda ({a_0}_y+{a_1}_x)+3 a_1\Lambda_x=0, \eqno(4.2)
$$
$$
2 (1+a_0+a_4)\Lambda_x+2 \Lambda ({a_1}_y+{a_4}_x)+3 a_1\Lambda_y=0, \eqno(4.3)
$$
$$
2\Lambda {a_4}_y + a_1\Lambda_x + 4 a_4 \Lambda_y=0. \eqno(4.4)
$$
For brevity we will write only some of the other equations:
$$
f_x = g_y, \eqno(4.5)
$$
$$
4 (c_0 - c_2)_x - 4 {c_1}_y -3 g (f_y+g_x) = 0, \eqno (4.6)
$$
$$
4 (c_0 - c_2)_y + 4 {c_1}_x +3 f (f_y+g_x) = 0, \eqno (4.7)
$$

$$
2 {d_0}_x - 2 {d_1}_y - c_1 (f_y+g_x) = 0, \eqno(4.8)
$$
$$
2 {d_0}_y + 2 {d_1}_x + (c_0 - c_2) (f_y+g_x) = 0, \eqno(4.9)
$$
$$
(d_0\Lambda)_x+(d_1\Lambda)_y=0, \eqno(4.10)
$$

$$
4 {e_0}_x - d_1 (f_y+g_x) = 0, \eqno(4.11)
$$
$$
4 {e_0}_y + d_0 (f_y+g_x) = 0. \eqno(4.12)
$$
Similar to the previous section, let us transform (4.6), (4.7) via (4.5) obtaining
$$
4 (c_0 - c_2) + \frac{3}{2}(f^2 - g^2)  = K_1, \eqno(4.13)
$$
$$
4 c_1 + 3 f g = K_2, \eqno(4.14)
$$
It follows from (4.8), (4.9) that
$$
2 \triangle d_0 + \left((c_0-c_2) (f_y+g_x)\right)_y - \left(c_1 (f_y+g_x)\right)_x = 0,
$$
$$
2 \triangle d_1 + \left((c_0-c_2) (f_y+g_x)\right)_x + \left(c_1 (f_y+g_x)\right)_y = 0,
$$
and due to (4.5), (4.13), (4.14) we obtain
$$
16 d_0 + 2 K_1 f - 2 K_2 g - f^3 +3 f g^2 = K_3. \eqno(4.15)
$$
$$
16 d_1 + 2 K_1 g + 2 K_2 f + g^3 -3 f^2 g = K_4, \eqno(4.16)
$$
Finally, it follows from (4.11), (4.12) that
$$
\left(d_0 (f_y+g_x)\right)_x + \left(d_1 (f_y+g_x)\right)_y = 0. \eqno(4.17)
$$
Similar to the previous section, due to (4.5), (4.15), (4.16) we verify that (4.17) is equivalent to
$$
fg (f^2 - g^2) - 2 K_1 fg + K_2 (f^2 - g^2) + K_3 g + K_4 f  = K_5. \eqno(4.18)
$$
Now we may rewrite (4.5) in the following way:
$$
g_y - f'(g) g_x = 0, \eqno(4.19)
$$
here $f'(g)$ can be found from (4.18).
By the same reasons as in the previous section the equation (4.19) does not admit non-constant global analytic periodic solutions.
Hence, at least one of the functions $f(x,y),\  g(x,y)$ is constant.
Without loss of generality we may assume that $g(x,y) = G$ is a constant.
Due to (4.5) we have $f(x,y) = f_1(y),$ hence the magnetic field $\Omega$ depends only on one variable:
$$\Omega (y) = \frac{1}{4}f_1'(y).$$
The relations (4.13), (4.14) take the form
$$
c_0 (x,y) - c_2(x,y) = \frac{1}{4} \left( K_1 + \frac{3}{2} G^2 - \frac{3}{2} f_1(y)^2 \right), \quad
c_1 (y) = \frac{1}{4} \left( K_2 - 3 G f_1(y) \right).
$$
It follows from (4.15), (4.16) that $d_0, \ d_1$ do not depend on $x.$ Due to (4.15), (4.16) we obtain that (4.11), (4.12) are equivalent to
$$64 {e_0}_x+\left( G^3 + 2 G K_1 - K_4 + 2 K_2  f_1(y) - 3 G f_1(y)^2 \right)f_1'(y) = 0, \eqno(4.20)$$
$$64 {e_0}_y+\left( 2 G^2 K_2 + K_3 - (3 G^2 + 2 K_1) f_1(y) + f_1(y)^3 \right)f_1'(y) = 0. \eqno(4.21)$$
It follows from (4.20) that $e_0(x,y) = e_0(y).$ We suppose $f_1(y)$ to be non-constant (since otherwise the magnetic field vanishes identically). Then (4.18) implies that $$G = K_2 = K_4 = K_5 = 0.$$
Now (4.21) gives
$$e_0(y) = -\frac{1}{64} K_3 f_1(y) +\frac{1}{64} K_1 f_1(y)^2 - \frac{1}{256}f_1(y)^4.$$
As in the previous section, (4.10) implies that $$\Lambda(x,y) = \Lambda(y).$$
Via arguments similar to those in the previous section we finally find all the others coefficients of the first integral $F_4:$
$$
a_0(y) = \frac{s_3+s_2\Lambda(y)-\Lambda(y)^2}{\Lambda(y)^2}, \quad a_1 = 0, \quad a_4(y) = 1 + a_0(y) - \frac{s_2}{\Lambda(y)},
$$
$$
b_2(y) = \frac{2 s_5 - s_2 f_1(y)}{2 \Lambda(y)}, \quad b_3(y) = 0, \quad c_2(y) = \frac{16 s_6 - 4 s_5 f_1(y) + s_2 f_1(y)^2}{16 \Lambda(y)},
$$
here $s_2, s_3, s_5, s_6$ are arbitrary constants.

Recall that the magnetic field has the form $\Omega(y) = \frac{f_1'(y)}{4}.$
There exists a linear in momenta first integral $F_1 = p_1 - \frac{f_1(y)}{4}$ on all energy levels.
One can check by straightforward calculations that
$$F_4 = -F_1^4 + 4 s_3 H^2 + 2 s_2 H F_1^2 + 2 s_5 H F_1 + 2 s_6 H + \frac{K_1}{4} F_1^2 + \frac{K_3}{16} F_1,$$
and it means that $F_4$ can be expressed in terms of integrals of lesser degree.
Due to remark in the section 2, this proves Theorem 2.

\section{General case}

Let us turn to the proof of the Theorem 3.
Suppose that there exists an additional polynomial in momenta first integral of degree $N$ on all energy levels:
$$F=\sum_{s=0}^{N} \sum_{k=0}^{s} a_{s,k}(x,y)p_1^{s-k}p_2^k. \eqno(5.1)$$
We assume $a_{s,k}(x,y)$ to be analytic periodic functions in both variables. Let us denote $$\alpha_j = a_{N-j,0}-a_{N-j,2}+a_{N-j,4- \ldots}, \qquad \beta_j = a_{N-j,1}-a_{N-j,3}+a_{N-j,5- \ldots},$$
here $a_{m,n}=0$ while $m < n.$
Notice that $\beta_N = 0$ and $$\alpha_0 = a_{N,0}-a_{N,2}+a_{N,4- \ldots}, \qquad \beta_0 = a_{N,1}-a_{N,3}+a_{N,5- \ldots}$$
are Kolokol'tsov constants which can be chosen in the following way:
$$
\alpha_0 = -1, \ \beta_0 = 0.
$$
Let us also denote $$\alpha_1 = f, \ \beta_1 = g.$$

\begin{lemma}
If the geodesic flow (1.1) on the 2-torus in a non-zero magnetic field admits the first integral (5.1) on all energy levels, then the following sequence of equations holds:
$$
N{\alpha_j}_x-N{\beta_j}_y-(N+1-j) \beta_{j-1} (f_y+g_x) = 0, \eqno(5.2)
$$
$$
N{\alpha_j}_y+N{\beta_j}_x+(N+1-j) \alpha_{j-1} (f_y+g_x) = 0, \eqno(5.3)
$$
here $j = 1, \ldots, N.$
Moreover, the following relations hold true: $$f_x = g_y, \qquad \Omega (x, y) = \frac{1}{N} (f_y+g_x), \qquad (\alpha_{N-1} \Lambda)_x+(\beta_{N-1} \Lambda)_y=0. \eqno(5.4)$$
\end{lemma}

Let us prove Lemma 3. The relation $\dot{F} = \{F, H\}_{mg} = 0,$ where $H = \frac{p_1^2+p_2^2}{2 \Lambda(x,y)}$ gives
$$
\frac{p_1^2}{2 \Lambda^2} (\Lambda_x F_{p_1} + \Lambda_y F_{p_2}) + \frac{p_2^2}{2 \Lambda^2} (\Lambda_x F_{p_1} + \Lambda_y F_{p_2}) + \frac{p_1}{\Lambda} F_x + \frac{p_2}{\Lambda} F_y + \frac{p_2}{\Lambda} \Omega F_{p_1} - \frac{p_1}{\Lambda} \Omega F_{p_2} = 0. \eqno(5.5)
$$
This is a polynomial in $p_1, p_2$ of degree $N+1$ which can be written in the following form
$$\{F, H\}_{mg} = \sum_{s=0}^{N+1} \sum_{k=0}^{s} M_{s,k}(x,y)p_1^{s-k}p_2^k = 0,$$
all the coefficients $M_{s,k}(x,y) \equiv 0,$ they can be found explicitly from (5.5):
$$
M_{s,k} = (s-k-1) a_{s-1,k} \frac{\Lambda_x}{2 \Lambda^2} + (k+1) a_{s-1,k+1} \frac{\Lambda_y}{2 \Lambda^2} + (s-k+1) a_{s-1,k-2} \frac{\Lambda_x}{2 \Lambda^2} +
$$
$$
(k-1) a_{s-1,k-1} \frac{\Lambda_y}{2 \Lambda^2} + \frac{(a_{s-1,k})_x}{\Lambda} +\frac{(a_{s-1,k-1})_y}{\Lambda} + (s-k+1) a_{s,k-1} \frac{\Omega}{\Lambda} - (k+1) a_{s,k+1} \frac{\Omega}{\Lambda} = 0.
$$
Here $a_{m,n}=0$ while $m < n, m<0, n < 0, m > N$ or $n > N.$
Introduce new functions $$P_{s,k} = \left((s-k-1) a_{s-1,k} + (s-k+1) a_{s-1,k-2}\right) \frac{\Lambda_x}{2 \Lambda^2},$$
$$Q_{s,k} = \left((k+1) a_{s-1,k+1} + (k-1) a_{s-1,k-1}\right) \frac{\Lambda_y}{2 \Lambda^2},$$ $$R_{s,k} = \frac{(a_{s-1,k})_x+(a_{s-1,k-1})_y}{\Lambda}, \quad T_{s,k} = \left((s-k+1) a_{s,k-1} - (k+1) a_{s,k+1}\right) \frac{\Omega}{\Lambda}.$$
Then $$M_{s,k}=P_{s,k} + Q_{s,k} + R_{s,k} + T_{s,k}.$$
Let us introduce $$U_s = M_{s,0} - M_{s,2} + M_{s,4} - \ldots, \qquad V_s = M_{s,1} - M_{s,3} + M_{s,5} - \ldots.$$
Notice that $$P_{s,0}-P_{s,2}+P_{s,4}- \ldots = 0, \qquad Q_{s,0}-Q_{s,2}+Q_{s,4}- \ldots = 0, $$
$$P_{s,1}-P_{s,3}+P_{s,5}- \ldots = 0, \qquad Q_{s,1}-Q_{s,3}+Q_{s,5}- \ldots = 0,$$
$$R_{s,0}-R_{s,2}+R_{s,4}- \ldots = \frac{1}{\Lambda} \left( (\alpha_{N-s+1})_x - (\beta_{N-s+1})_y \right),$$
$$R_{s,1}-R_{s,3}+R_{s,5}- \ldots = \frac{1}{\Lambda} \left( (\alpha_{N-s+1})_y + (\beta_{N-s+1})_x \right),$$
$$T_{s,0}-T_{s,2}+T_{s,4}- \ldots = -s \frac{\Omega}{\Lambda} \beta_{N-s}, \quad T_{s,1}-T_{s,3}+T_{s,5}- \ldots = s \frac{\Omega}{\Lambda} \alpha_{N-s}.$$
Put $s = N+1-j,$ then the relations $U_s = 0, V_s = 0$ give
$$
{\alpha_j}_x-{\beta_j}_y-(N+1-j) \Omega \beta_{j-1} = 0,
$$
$$
{\alpha_j}_y+{\beta_j}_x+(N+1-j) \Omega \alpha_{j-1} = 0.
$$
Since $\alpha_0 = -1, \beta_0 = 0, \alpha_1=f, \beta_1=g,$ for $j=1$ we obtain $f_x = g_y, \Omega = \frac{1}{N} (f_y+g_x).$

It is left to notice that $M_{2,0}+M_{2,2} = 0$ implies $(\alpha_{N-1} \Lambda)_x+(\beta_{N-1} \Lambda)_y=0.$ Lemma 3 is proved.

By cross differentiation of (5.2), (5.3) we obtain:
$$
N \triangle \alpha_j + (N+1-j) \left( (\alpha_{j-1}(f_y+g_x))_y - (\beta_{j-1}(f_y+g_x))_x  \right) = 0. \eqno(5.6)
$$
$$
N \triangle \beta_j + (N+1-j) \left( (\alpha_{j-1}(f_y+g_x))_x + (\beta_{j-1}(f_y+g_x))_y  \right) = 0. \eqno(5.7)
$$
Now we claim that, due to (5.6), (5.7), $\alpha_j, \beta_j$ are polynomials in $f,g$ with constant coefficients for any $j=1, \ldots, N.$ We shall prove this statement in the appendix (see the Theorem 4). Recall that $\beta_N \equiv 0.$ Taking this into account we may conclude that there is a nontrivial polynomial $P$ such that $P(f,g)=0.$
Due to $f_x = g_y$ we obtain:
$$P'_f g_y+P'_g g_x=0, \qquad P'_f f_y+P'_g f_x=0.$$
This implies that characteristics of $f$ and $g$ coincide. By similar reasons as in section 3 it can be shown that one of the functions $f, g$ must be constant. Assume that $g(x,y) \equiv G$ where $G$ is a constant. Then, due to (5.4) we have $$f(x,y) = f(y). \qquad \Omega(x,y) = \frac{f'(y)}{N}$$
and it means that the magnetic field depends only on one variable $y.$
Notice that $\alpha_j = \alpha_j(f,g), \beta_j = \beta_j(f,g).$ It means that
$\alpha_j = \alpha_j(y), \beta_j = \beta_j(y), j = 0, \ldots, N.$ Due to this, at $j=N$ in (5.2) we have
$$
\beta_{N-1}(y) f_y = 0. \eqno(5.8)
$$
We recall that $\beta_{N-1}(y)$ is a polynomial in $f(y).$ Let us write it down in the following way: $\beta_{N-1}(y) = \omega_k f^k(y) + \omega_{k-1} f^{k-1}(y) + \ldots + \omega_1 f(y) +\omega_0,$ where $\omega_j$ are constants. Then (5.8) is equivalent to $$\frac{\omega_k}{k+1} f^{k+1}(y) + \ldots + \frac{\omega_1}{2} f^2(y) + \omega_0 f(y) = c_0$$ with an arbitrary constant $c_0.$ Notice that $f(y)$ is not a constant (since otherwise the magnetic field vanishes identically). Due to this, the last equation implies that $\omega_j = 0$ for any $j=0, \ldots, k.$
Thus we have $\beta_{N-1}(y) \equiv 0.$
At $j=N-1$ in (5.2) we have
$$
\beta_{N-2}(y) f_y = 0.
$$
from which it follows that $\beta_{N-2}(y) \equiv 0.$
Proceeding in this way we obtain $$\beta_N = \beta_{N-1} = \ldots = \beta_2 = \beta_1 = g \equiv G = 0.$$
Then, (5.3) implies that
$$
N{\alpha_j}_y+(N+1-j) \alpha_{j-1} f_y = 0. \eqno(5.9)
$$
We recall that $\alpha_j$ is a polynomial in $f$ for any $j = 0, \ldots, N.$ Suppose that there exist $k, 1 < k \leq N$ such that $\alpha_k(y) \equiv 0.$ Then, due to (5.9) we have $\alpha_{k-1}(y) \equiv 0.$ Proceeding in this way we obtain $\alpha_1(y) = f(y) \equiv 0$ but in this case the magnetic field vanishes identically. Thus we have $\alpha_{j} \neq 0$ for any $j.$

Finally, due (5.4), the following equation holds true (compare with (3.17), (4.10)):
$$(\alpha_{N-1} \Lambda)_x+(\beta_{N-1} \Lambda)_y=0.$$
Since $\beta_{N-1} \equiv 0$ and $\alpha_{N-1}(x,y) = \alpha_{N-1}(y) \neq 0,$ we have
$$\alpha_{N-1}(y) \Lambda_x=0. \eqno(5.10)$$
Suppose that there exists a point $q=(x_0,y_0)$ such that $\Lambda_x \neq 0$ at $q.$ It means that there is a neighborhood $U_q$ of $q$ such that $\Lambda_x \neq 0$ in $U_q.$ Then, due to (5.10) $\alpha_{N-1}(y) \equiv 0$ in $U_q.$ Recall that $f(y)$ is an analytic function and $\alpha_{N-1}(y)$ is a polynomial in $f(y)$ with constant coefficients. Due to this, if $\alpha_{N-1}(y)$ is a constant on an open subset, then it is identically constant on the whole domain: $\alpha_{N-1}(y) \equiv 0.$ But this case has already been considered: namely, in this case $f(y) \equiv const$ and $\Omega(y) \equiv 0.$ So we have $\Lambda_x \equiv 0$ from what follows that the metrics depends only on one variable:
$$\Lambda(x,y) = \Lambda(y).$$
Hence we obtain that the metrics $\Lambda$ and the magnetic field $\Omega$ are functions of only one variable. But in this case the linear in momenta first integral $F_1 = p_1-\frac{f(y)}{N}$ exists on all energy levels. Due to lemma 2, this proves the Theorem 3.

\section{Appendix}
In this section we prove that $\alpha_j, \beta_j$ depend polynomially on $f,g$ for any $j=1, \ldots, N.$ Recall that $\alpha_0=-1$, $\beta_0=0$, $\alpha_1=f$ and $\beta_1=g$. Moreover, functions $\alpha_j$ and $\beta_j$ satisfy the following recurrence relations:
$$
N\triangle\alpha_j+(N+1-j)([\alpha_{j-1}(f_y+g_x)]_y-[\beta_{j-1}(f_y+g_x)]_x)=0, \eqno(6.1)
$$
$$
N\triangle\beta_j+(N+1-j)([\alpha_{j-1}(f_y+g_x)]_x+[\beta_{j-1}(f_y+g_x)]_y)=0. \eqno(6.2)
$$
We also have $f_x=g_y.$

Let us define the polynomials $A_n$ and $B_n$ for $n\geq 0$ as follows:
$$A_n=\frac{1}{(n+1)!} \sum_{k=0}^{\lfloor\frac{n+1}{2}\rfloor} f^{n-2k+1}g^{2k}(-1)^k\binom{n+1}{2k},$$
$$B_n=\frac{1}{(n+1)!} \sum_{k=1}^{\lfloor\frac{n+2}{2}\rfloor} f^{n-2k+2}g^{2k-1}(-1)^{k+1}\binom{n+1}{2k-1}.$$
For $n=0$ and $n=1$ we have $A_0 = f, B_0 = g, A_1 = \frac{f^2-g^2}{2}, B_1 = fg.$
Notice that
$$
A_n=\frac{1}{2(n+1)!}[(f+ig)^{n+1}+(f-ig)^{n+1}], \eqno(6.3)
$$
$$
B_n=\frac{1}{2i(n+1)!}[(f+ig)^{n+1}-(f-ig)^{n+1}]. \eqno(6.4)
$$

\begin{lemma}
Polynomials $A_n$ and $B_n$ satisfy the following recurrence relations:
$$
A_n(f_y+g_x)=(A_{n+1})_y+(B_{n+1})_x, \eqno(6.5)
$$
$$
B_n(f_y+g_x)=-(A_{n+1})_x+(B_{n+1})_y. \eqno(6.6)
$$
\end{lemma}

Let us prove lemma 4. Since (6.5) and (6.6) can be proved analogously, we shall prove only (6.5).

Using the equalities (6.3), (6.4), we obtain
$$(A_{n+1})_y=\frac{1}{2(n+1)!}[(f+ig)^{n+1}(f_y+ig_y)+(f-ig)^{n+1}(f_y-ig_y)],$$
$$(B_{n+1})_x=\frac{1}{2i(n+1)!}[(f+ig)^{n+1}(f_x+ig_x)-(f-ig)^{n+1}(f_x-ig_x)].$$
Hence we have
$$(A_{n+1})_y=A_n\cdot f_y-B_n \cdot g_y, \qquad
(B_{n+1})_x=B_n\cdot f_x+A_n \cdot g_x.$$
Since $f_x=g_y$, we obtain $(A_{n+1})_y+(B_{n+1})_x=A_n(f_y+g_x).$ Lemma 4 is proved.

\begin{proposition}
The following relations hold true:
$$\alpha_2=\frac{1}{N}a_2-\frac{N-1}{N}A_1, \quad \beta_2=\frac{1}{N}b_2-\frac{N-1}{N}B_1.$$
\end{proposition}

Notice that $\alpha_1=A_0, \ \beta_1=B_0.$ Hence $\alpha_1(f_y+g_x)=A_0(f_y+g_x), \ \beta_1(f_y+g_x)=B_0(f_y+g_x).$ Then, due to (6.5), (6.6) we obtain $$[\alpha_1(f_y+g_x)]_y-[\beta_1(f_y+g_x)]_x=\triangle A_1.$$
Then $\alpha_2$ satisfies $$N\triangle\alpha_2+(N-1)\triangle A_1=0.$$
Hence, $N\alpha_2+(N-1)A_1=a_2,$ where $a_2$ is a constant. So we obtain $\alpha_2=\frac{1}{N}a_2-\frac{N-1}{N}A_1$. Function $\beta_2$ can be found analogously.

\begin{proposition}
The following relations hold true:
$$\alpha_3=\frac{1}{N}a_3-\frac{N-2}{N^2}(a_2A_0-b_2B_0)+\frac{(N-2)(N-1)}{N^2}A_2,$$ $$\beta_3=\frac{1}{N}b_3-\frac{N-2}{N^2}(a_2B_0+b_2A_0)+\frac{(N-2)(N-1)}{N^2}B_2.$$
\end{proposition}
Due to Proposition 1 we have $$\alpha_2(f_y+g_x)=\frac{1}{N}a_2(f_y+g_x)-\frac{N-1}{N}A_1(f_y+g_x),$$
$$\beta_2(f_y+g_x)=\frac{1}{N}b_2(f_y+g_x)-\frac{N-1}{N}B_1(f_y+g_x).$$
Then (6.5), (6.6) and $f_x=g_y$ imply that
$$[\alpha_2(f_y+g_x)]_y-[\beta_2(f_y+g_x)]_x=\frac{1}{N}(a_2\triangle f-b_2 \triangle g)-\frac{N-1}{N}\triangle A_2.$$
Then $\alpha_3$ satisfies
$$N\triangle\alpha_3+(N-2)\triangle(\frac{1}{N}(a_2f-b_2 g)-\frac{N-1}{N} A_2)=0.$$
Hence, $$N\alpha_3+\frac{N-2}{N}(a_2f-b_2 g)-\frac{(N-2)(N-1)}{N} A_2=a_3,$$ where $a_3$ is a constant. So we obtain
$$\alpha_3=\frac{1}{N}a_3-\frac{N-2}{N^2}(a_2A_0-b_2B_0)+\frac{(N-2)(N-1)}{N^2}A_2.$$ Function $\beta_3$ can be found analogously.

Introduce $c_j=(-1)^j\frac{(N-j)\ldots(N-1)}{N^j}$.
\begin{theorem}
For $j\geq3$ functions $\alpha_j, \beta_j$ depend polynomially on $f,g$ and can be found in the following way:
$$
\alpha_{j}=\frac{1}{N}a_j+\sum_{i=0}^{j-3}(m_{i}(j)A_i-n_i(j)B_i)+c_{j-1}A_{j-1}, \eqno(6.7)
$$
$$
\beta_{j}=\frac{1}{N}b_j+\sum_{i=0}^{j-3}(m_{i}(j)B_i+n_i(j)A_i)+c_{j-1}B_{j-1}, \eqno(6.8)
$$
where $a_j, b_j$ are arbitrary constants and $m_i(j), n_i(j) \in \mathbb{R}.$
\end{theorem}

Let us prove this theorem by induction on $j.$ The base for $j=3$ is proved in Proposition 2. The proof is analogous for $\alpha_j$ and $\beta_j$ hence let us prove the theorem only for $\alpha_j.$ Prove the induction step.

It follows from (6.5) and (6.6) that
$$
[A_i(f_y+g_x)]_y-[B_i(f_y+g_x)]_x=\triangle A_{i+1}, \eqno(6.9)
$$
$$
[B_i(f_y+g_x)]_y+[A_i(f_y+g_x)]_x=\triangle B_{i+1}. \eqno(6.10)
$$

By induction we have $$\alpha_j(f_y+g_x)=\frac{1}{N}a_j(f_y+g_x)+\sum_{i=0}^{j-3}(m_{i}(j)A_i(f_y+g_x)-n_i(j)B_i(f_y+g_x))+c_{j-1}A_{j-1}(f_y+g_x),$$
$$\beta_j(f_y+g_x)=\frac{1}{N}b_j(f_y+g_x)+\sum_{i=0}^{j-3}(m_{i}(j)B_i(f_y+g_x)+n_i(j)A_i(f_y+g_x))+c_{j-1}B_{j-1}(f_y+g_x).$$
Hence, due to (6.9), (6.10) and $f_x=g_y$ we obtain $$[\alpha_j(f_y+g_x)]_y-[\beta_j(f_y+g_x)]_x=\frac{1}{N}(a_j\triangle f-b_j\triangle g)+\sum_{i=0}^{j-3}(m_{i}(j)\triangle A_{i+1}-n_i(j)\triangle B_{i+1})+c_{j-1}\triangle A_{j}.$$
So $$[\alpha_j(f_y+g_x)]_y-[\beta_j(f_y+g_x)]_x=\sum_{i=0}^{j-2}(m'_{i}\triangle A_{i}-n'_i\triangle B_{i})+c_{j-1}\triangle A_{j},$$ where $m'_i$ and $n'_i$ are constants.
Then the recurrence relation for $\alpha_{j+1}$ has the form $$N\triangle \alpha_{j+1}+(N-j)\triangle (\sum_{i=0}^{j-2}(m'_{i}A_{i}-n'_i B_{i})+c_{j-1}A_{j})=0.$$
Hence, $$N\alpha_{j+1}+(N-j)\sum_{i=0}^{j-2}(m'_{i}A_{i}-n'_i B_{i})+(N-j)c_{j-1}A_{j}=a_{j+1},$$ where $a_{j+1}$ is an arbitrary constant.
So we obtain $$\alpha_{j+1}=\frac{1}{N}a_{j+1}-\frac{N-j}{N}\sum_{i=0}^{j-2}(m'_{i}A_{i}-n'_i B_{i})-\frac{N-j}{N}c_{j-1}A_{j}.$$
The induction step is proved and this completes the proof of the Theorem 4.

\

S.V. Agapov

Sobolev Institute of Mathematics, Novosibirsk, Russia,


Novosibirsk State University, Novosibirsk, Russia

agapov@math.nsc.ru, \ \ agapov.sergey.v@gmail.com

\vspace{5mm}

A.A. Valyuzhenich

Sobolev Institute of Mathematics, Novosibirsk, Russia

graphkiper@mail.ru

\end{document}